\documentclass[12pt,a4paper]{article}

\usepackage{amsmath}
\usepackage{amsfonts}
\usepackage{latexsym}
\newcommand{\proof}{\noindent{\it Proof: }}
\newcommand{\proofbox}{\mbox{ $\Box$}\\}

\newcommand{\R}{\mathbb{R}}
\newcommand{\Z}{\mathbb{Z}}

\begin{document}

\author{Franck Barthe, K\'aroly J. B\"or\"oczky\footnote{Supported by
 the FP7 IEF grant GEOSUMSETS}, Matthieu Fradelizi }

\title{Stability of the functional forms of the Blaschke-Santal\'o
inequality}

\newtheorem{lemma}{Lemma}[section]
\newtheorem{theo}[lemma]{Theorem}
\newenvironment{thm}[3][]{\pagebreak[3] \medskip {\bf \noindent Theorem}
                               {#1} {#2}{\it #3}
                             }{}
\newtheorem{defi}[lemma]{Definition}
\newtheorem{coro}[lemma]{Corollary}
\newtheorem{conj}[lemma]{Conjecture}
\newtheorem{prop}[lemma]{Proposition}
\newtheorem{remark}[lemma]{Remark}
\newtheorem{example}[lemma]{Example}

\maketitle


\begin{abstract}
Stability versions
of the functional forms of the Blaschke-Santal\'o inequality
due to Ball, Artstein-Klartag-Milman, Fradelizi-Meyer and Lehec
are proved.
\end{abstract}

\newpage

\section{Introduction}

For general references about convex bodies,
see  P.M.~Gruber \cite{Gru07} or R.~Schneider \cite{Sch93},
and for a survey on related geometric inequalities,
see E.~Lutwak \cite{Lut93}.
We write $0$ to denote the origin of $\R^n$,
$\langle \cdot,\cdot\rangle$ to denote
the standard scalar product, $|\cdot|$
to denote the corresponding $l_2$-norm, and $V(\cdot)$ to denote
volume (Lebesgue-measure). Let $B^n$ be the unit Euclidean ball, and let $S^{n-1}=\partial B^n$.
A convex body $K$ in $\R^n$ is a compact convex set
with non--empty interior.
If $z\in{\rm int}K$,
then the polar of $K$ with respect to $z$ is the convex body
$$
K^z=\{x\in\R^n:\,\langle x-z,y-z\rangle\leq 1\mbox{ for any $y\in K$}\}.
$$
From Hahn-Banach's theorem in $\R^n$,  $(K^z)^z=K$.
According to L.A.~Santal\'o \cite{San49}
 (see also M.~Meyer and A.~Pajor \cite{MeP90}), there exists a unique
$z\in{\rm int}K$ minimizing the volume product
$V(K)V(K^z)$, which is called the Santal\'o point
of $K$. In this case $z$ is the centroid of $K^z$.
The Blaschke-Santal\'o
inequality states that if $z$ is the  Santal\'o point
(or centroid) of $K$, then
\begin{equation}
\label{BS}
V(K)V(K^z)\leq V(B^n)^2,
\end{equation}
with equality if and only if $K$ is an ellipsoid.
The inequality was proved by W.~Blaschke \cite{Bla17} (available also in \cite{Bla22}) for $n\leq 3$,
 and by L.A.~Santal\'o \cite{San49} for all $n$.
The case of equality was characterized by
J.~Saint-Raymond \cite{Sai81} among $o$-symmetric convex bodies,
and by C.M.~Petty \cite{Pet85} among all convex bodies
(see also D.~Hug \cite{Hug96}, E.~Lutwak \cite{Lut91},
M.~Meyer and A.~Pajor \cite{MeP90}, and M.~Meyer and S.~Reisner \cite{MeR06}
 for simpler proofs).

To state functional versions of the Blaschke-Santal\'o inequality, let us first recall that the usual definition of the Legendre transform
of a function $\varphi:\R^n\to \R\cup\{+\infty\}$ at
$z\in\R^n$ is defined by
$$
{\cal L}_z\varphi(y)=\sup_{x\in\R^n}\{\langle x-z,y-z\rangle-\varphi(x)\},
\mbox{ \ for $y\in\R^n$}
$$
and that the function ${\cal L}_z\varphi:\R^n\to \R\cup\{+\infty\}$ is always convex and lower semicontinuous.
If $\varphi$ is convex, lower semicontinuous and $\varphi(z)<+\infty$ then ${\cal L}_z{\cal L}_z\varphi=\varphi$.\\

Subsequent work by K.M.~Ball \cite{Bal},
S.~Artstein-Avidan, B.~Klartag, V.D.~Milman \cite{AKM04},
M.~Fradelizi, M.~Meyer \cite{FrM07} and J.~Lehec \cite{Leh1,Leh2}  lead to the functional
version of the Blaschke-Santal\'o inequality 
(see \cite{Bal} and \cite{AKM04} for the relation
between the functional version and the original Blaschke-Santal\'o inequality).

\begin{thm}[\cite{Bal, AKM04, FrM07, Leh1, Leh2}]{}{
Let $\varrho:\R\to \R_+$ be a log-concave non-increasing function and 
 $\varphi:\R^n\to \R$ be measurable then
$$
\inf_{z\in\R^n}\int_{\R^n}\varrho(\varphi(x))\,dx
\int_{\R^n}\varrho({\cal L}_z\varphi(x))\,dx\leq
\left(\int_{\R^n}\varrho(|x|^2/2)\,dx\right)^2.
$$
If $\varrho$ is decreasing there is equality if and only if there exist
$a,b,c\in\R$, $a<0$, $z\in\R^n$ and a positive definite matrix $T:\R^n\to\R^n$, such that
$$
\varphi(x)=\frac{|T(x+z)|^2}2+c
\mbox{ \ for $x\in\R^n$},
$$
and moreover either $c=0$, or $\varrho(t)=e^{at+b}$ for $t>-|c|$.
}\end{thm}

\vspace{1em}
Here we prove a stability version of this inequality.

\begin{theo}
\label{Legendrestab}
 Let  $\varrho:\R\to \R_+$ be a  log-concave and decreasing function with
$\int_{\R_+}\varrho<+\infty$.
Let $\varphi:\R^n\to \R$ be measurable. Assume that   
for some  $\varepsilon\in(0,\varepsilon_0)$ and for all $z\in\R^n$, the following inequality holds:
$$
\int_{\R^n}\varrho(\varphi(x))\,dx \int_{\R^n}
\varrho({\cal L}_z\varphi(x))\,dx>(1-\varepsilon)
\left(\int_{\R^n}\varrho(|x|^2/2)\,dx\right)^2
$$
\begin{enumerate}
\item If $\varphi$ is convex, then 
there exist some $z\in\R^n$,
$c\in\R$  and a positive definite matrix $T:\R^n\to\R^n$,
such that
$$
\int_{R(\varepsilon)B^n}
\left|\frac{|x|^2}2+c-\varphi(Tx+z)\right|\,dx
<\eta\varepsilon^{\frac1{129n^2}},
$$
where $\lim_{\varepsilon\to 0}R(\varepsilon)=+\infty$, and
$\varepsilon_0,\eta,R(\varepsilon)$ depend on $n$ and $\varrho$.
\item If $\varphi$ is only assumed to be measurable then a weaker version holds: There exists $z,c,T$ as above 
and $\Psi\subset R(\varepsilon)B^n$
 such that
$$
\int_{R(\varepsilon)B^n\backslash \Psi}
\left|\frac{|x|^2}2+c-\varphi(Tx+z)\right|\,dx
<\eta\varepsilon^{\frac1{129n^2}},
$$
and $V(\Psi\cap R B^n)\leq \eta\sqrt{\varepsilon}\,R^n$
for any $R\in[R_0,R(\varepsilon)]$, where $R_0>0$ depends only on $\varrho$.
\end{enumerate}
\end{theo}

\begin{remark}
One cannot expect the $L_1$-distance between $\varphi$ and $\frac{|T(x+z)|^2}2+c$
to be small on the whole $\R^n$. For instance, if $\varrho(t)=e^{-t}$,
and for small $\varepsilon>0$, $\varphi(x)=|x|^2/2$ if $|x|\leq |\log\varepsilon|$,
and $\varphi(x)=+\infty$ if $|x|> |\log\varepsilon|$, then, of course, for any $c$ and $T$ the function 
$x\mapsto \frac{|T(x+z)|^2}2+c-\varphi(x)$ is not in $L_1$, but 
$$
\int_{\R^n}\varrho(\varphi(x))\,dx \int_{\R^n}
\varrho({\cal L}_z\varphi(x))\,dx>(1-O(\varepsilon|\log\varepsilon|^{n-1}))
\left(\int_{\R^n}\varrho(|x|^2/2)\,dx\right)^2
$$
for all $z\in\R^n$.

In addition, if $\varphi$ is only assumed to be measurable, then we
may choose it to be infinity on a ball of small enough measure, and set 
$\varphi(x)=|x|^2/2$ on the complement.

On the other hand, most probably the exponent $\frac1{129n^2}$ in Theorem~\ref{Legendrestab} can be
exchanged into some positive absolute constant.
\end{remark}

As a matter of fact, the above functional form of the Blaschke-Santal\'o inequality
deduces from the following more general inequality, which is the result of different contributions as explained below

\begin{thm}[\cite{Bal, AKM04, FrM07, Leh1, Leh2}]{}{
\label{FraMey}
For any measurable $f:\R^n\to \R_+$ with positive integral
there exists a particular point $z\in\R^n$ attached to $f$ such that
 if measurable functions $\varrho:\R_+\to \R_+$
and $g:\R^n\to \R_+$ with positive integrals
satisfy
$$
f(x)g(y)\leq\varrho^2(\langle x-z,y-z\rangle),
$$
for every $x,y\in\R^n$ with $\langle x-z,y-z\rangle>0$, then
$$
\int_{\R^n}f(x)\,dx \int_{\R^n} g(x)\,dx\leq
\left(\int_{\R^n}\varrho(|x|^2)\,dx\right)^2.
$$
Equality holds for this $z$ if and only if there exist
$\tilde{\varrho}:\R_+\to \R_+$, $\xi>0$ and a positive definite matrix
$T:\R^n\to\R^n$, such that $\tilde{\varrho}(e^t)$ is log-concave,
and for a.e. $x\in \R^n$ and $s\in\R_+$, we have
$$
\varrho(s)=\tilde{\varrho}(s),
\mbox{ \ }f(x)=\xi\,\tilde{\varrho}(|T(x-z)|^2)
\mbox{ \ \ and \ }g(x)=\xi^{-1}\tilde{\varrho}(|T^{-1}(x-z)|^2).
$$
}\end{thm}

K.M. Ball \cite{Bal} initiated the study of such inequalities, established the case of even functions $f$ and proved that, in this case, $z$ can be chosen to be the origin.
If $\varrho(t)=e^{-t}$, S. Artstein, B. Klartag, V.D. Milman \cite{AKM04}
showed that one can choose $z$ to be the mean of $f$ for any $f$. 
For any measurable $\varrho$ but for log-concave functions $f$, 
M. Fradelizi, M. Meyer \cite{FrM07} constructed the
suitable $z$ in the following way.
For any $z\in\R^n$, let
$$
K_{f,z}=\left\{x\in\R^n:\,\int_0^{+\infty} r^{n-1}f(z+rx)\,dx\geq 1\right\},
$$
which is convex according to K.M. Ball \cite{Bal88}.
M. Fradelizi and M. Meyer \cite{FrM07}
proved that there exists a $z\in\R^n$, such that
the centre of mass of $K_{f,z}$ is the origin and that
this $z$ works.
Finally  J. Lehec gave a direct and different proof of the general theorem in \cite{Leh2}. 
He established the existence of a so-called Yao-Yao center for any measurable $f$ 
and that this point $z$ works also. 

We also give a stability version of this more general form of the Blaschke-Santal\'o inequality.

\begin{theo}
\label{FraMeystab}
 If some log-concave functions $\varrho:\R_+\to \R_+$ and
 $f,g:\R^n\to \R_+$ with positive integrals satisfy that
$\varrho$ is non-increasing, the centre
of mass of $K_{f,z}$ is the origin for some $z\in\R^n$, and
$$
f(x)g(y)\leq\varrho^2(\langle x-z,y-z\rangle)
$$
for every $x,y\in\R^n$ with $\langle x-z,y-z\rangle>0$, if
 moreover for $\varepsilon>0$,
$$
(1+\varepsilon)\int_{\R^n}f(x)\,dx \int_{\R^n} g(x)\,dx\geq
\left(\int_{\R^n}\varrho(|x|^2)\,dx\right)^2,
$$
then there exist $\xi>0$ and a positive definite matrix
$T:\R^n\to\R^n$, such that
\begin{eqnarray*}
\int_{\R^n}\left|\varrho(|x|^2)-\xi\,f(Tx+z)\right|\,dx
&<&\gamma\varepsilon^{\frac1{32n^2}} \cdot
\int_{\R_+}r^{n-1}\varrho(r^2)\,dr\\
\int_{\R^n}\left|\varrho(|x|^2)-\xi^{-1}g(T^{-1}x+z)\right|\,dx
&<&\gamma\varepsilon^{\frac1{32n^2}}\cdot\int_{\R_+}r^{n-1}\varrho(r^2)\,dr ,
\end{eqnarray*}
where $\gamma$ depends only on $n$.
\end{theo}

We strongly believe that the power  $\frac1{32n^2}$
occurring in Theorem \ref{FraMeystab}
can be chosen to be a positive absolute constant.

In this note, the implied constant in $O(\cdot)$
depends only on the dimension $n$.

\section{Stability of the Borell and
the Blaschke-Santal\'o inequalities}

C. Borell \cite{Bor75} pointed out the following version
of the Pr\'ekopa-Leindler inequality:

\begin{theo}[Borell]
\label{PLB}
If $M,F,G:\R_+\to\R_+$ are integrable functions
with positive integrals, and   $M(\sqrt{rs})\geq \sqrt{F(r)G(s)}$ for
$r,s\in\R_+$,  then
$$
\int_{\R_+}F \cdot \int_{\R_+} G\leq \left(\int_{\R_+}  M\right)^2.
$$
\end{theo}

Recently the following stability estimate
has been obtained in K.M.~Ball, K.J.~B\"or\"oczky \cite{BB2}.
We note that if $M:\R_+\to \R_+$ is log-concave and non-increasing,
then $M(e^t)$ is log-concave on $\R$.

\begin{theo}[Ball, B\"or\"oczky]
\label{PLBstab}
There exists a positive absolute constant $c$ with
the following property:
If $M,F,G:\R_+\to\R_+$ are integrable functions
with positive integrals such that
$M(e^t)$ is log-concave,  $M(\sqrt{rs})\geq \sqrt{F(r)G(s)}$ for
$r,s\in \R_+$,  and
$$
\left(\int_{\R_+}  M\right)^2\leq (1+\varepsilon)
\int_{\R_+}F \cdot \int_{\R_+} G,
$$
for some $\varepsilon>0$, then
there exist $a,b>0$, such that
\begin{eqnarray*}
\int_{\R_+}|a\,F(bt)-M(t)|\,dt&\leq &
c\cdot\varepsilon^{\frac5{16}}\cdot
\int_{\R_+} M(t)\,dt \\
\int_{\R_+}|a^{-1}G(b^{-1}t)-M(t)|\,dt&\leq &
c\cdot\varepsilon^{\frac5{16}}\cdot
\int_{\R_+} M(t)\,dt.
\end{eqnarray*}
\end{theo}

For a stability version of the Blaschke-Santal\'o inequality,
we use the Banach-Mazur distance of two convex bodies $M$ and $K$,
which is defined by
$$
\delta_{\rm BM}(K,M) =\min\{ \ln\lambda:\,
K-x\subset \Phi(M)\subset \lambda(K-x)  \mbox{ for }
\Phi\in{\rm GL}(n),x\in\R^n\}.
$$
Improving on K.J. B\"or\"oczky \cite{Bor},
the paper \cite{BB2} also established the following.

\begin{theo}[Ball,B\"or\"oczky]
\label{BSstab}
If $K$ is a  convex body in $\R^n$, $n\geq 3$, with centroid $z$,
and
$$
V(K)V(K^z)>(1-\varepsilon)V(B^n)^2\mbox{ \ for some $\varepsilon\in(0,\frac12)$,}
$$
then for some $\gamma>0$ depending only on $n$, we have
$$
\delta_{\rm BM}(K,B^n)<\gamma\,\varepsilon^{\frac1{5n}}.
$$
\end{theo}

We note that according to K.M. Ball \cite{Bal}, Borell's inequality Theorem~\ref{PLB}
can be used to prove the Blaschke-Santal\'o inequality.
In particular, \cite{BB2} proves Theorem~\ref{BSstab}
via Theorem~\ref{PLBstab}.


\section{Proof of Theorem~\ref{FraMeystab}}

 Before proving Theorem~\ref{FraMeystab},
we verify first a simple property of log-concave
functions, then show that the centroid is a reasonable
centre for the Banach-Mazur distance from ellipsoids.

\begin{prop}
\label{logconcave}
If $h,\omega:\R\to \R_+$ are log-concave, $\omega$ is even, and
$$
\int_{\R}|r|^{n-1}|h(r)-\omega(r)|\,dr\leq
\varepsilon\int_\R |r|^{n-1}\omega(r)\,dr
$$
for some $\varepsilon\in (0,(250n)^{-(n+1)})$, then
$|h(0)-\omega(0)|\leq 250n\varepsilon^{\frac1{n+1}}\cdot \omega(0)$.
\end{prop}

\proof We may assume that $\omega(0)=1$ and
$\int_{\R} \omega(r)\,dr=1$, and hence $\omega(r)\leq 1$
for all $r$. First, we put forward a few useful facts about the function $\omega$.

 Following ideas from K.M. Ball and K.J. B\"or\"oczky \cite{BB1}, let us prove first that
there exists some $r_0\geq \frac{1}{2}$ such that
$\omega(r)\geq e^{-2|r|}$ if $|r|\leq r_0$, and $\omega(r)\leq e^{-2|r|}$
if $|r|\geq r_0$. 
For this,  notice that since $\int_{\R_+}\omega(r)dr=\frac{1}{2}=\int_{\R_+}e^{-2r}dr$ and  $\log\omega$ is concave there exists  $r_0>0$ satisfying the required property (and $r_0$ is unique unless $\omega(r)=e^{-2|r|}$ for all $r$. In this very specific case, we set arbitrarily
$r_0=1/2$).

 Now let us prove that $r_0\ge 1/2$.
We define $\omega^{-1}(t)=\sup\{r\ge 0; \omega(r)\ge t\}$. The hypotheses on $\omega$ imply that the support of $\omega^{-1}$ is $[0,1]$ and $\int_0^1\omega^{-1}(t)dt=1/2$. From Jensen's inequality one deduces that
$$
\omega\left(\frac{1}{2}\right)=\omega\left(\int_0^1\omega^{-1}(t)dt\right)\ge e^{\int_0^1\log(t)dt}=e^{-1}.
$$
Since $\omega(0)=1$, it follows from the log-concavity of $\omega$ that $\omega(r)\ge e^{-2r}$, if $|r|\leq 1/2$. This proves the claim.

 In particular, the latter exponential lower bound on $\omega$ implies that
\begin{equation}
\label{omegaest}
\omega(r)\geq 1-2|r| \quad
\mathrm{if} \quad |r|\leq\frac12.
\end{equation}
The fact that the graphs of $\omega$ and $r\mapsto e^{-2|r|}$ cross only once on $\R_+$ implies the following useful bound
\begin{equation}
\label{omegaint}
\int_{\R} r^{n-1}\omega(r)\,dr \leq 2\int_{\R_+} r^{n-1}e^{-2r}\,dr
=\frac{(n-1)!}{2^{n-1}}\leq n^{n+1}.
\end{equation}

Next, we study the function $h$.
Let $a_i=in\varepsilon^{\frac1{n+1}}$ for $i\in\Z$.
We claim that there exist two ind ices  $i\in\{1,\ldots,5\}$,
such that
\begin{equation}
\label{aiaj}
1-11n\varepsilon^{\frac1{n+1}}\leq
h(a_i) \leq 1+n\varepsilon^{\frac1{n+1}}.
\end{equation}
Suppose that (\ref{aiaj}) does not hold.
Since $h$ is non-decreasing and then non-increasing, there exists  $k\in\{1,2,3,4\}$
such that $h$ is monotone on $[a_k,a_{k+1}]$, and  $h(a_k)$ and $h(a_{k+1})$ are outside and on the same side
of the interval $[1-11n\varepsilon^{\frac1{n+1}}, 1+n\varepsilon^{\frac1{n+1}}]$.
Consequently, for this value of $k$,
 either $h(r)<1-11n\varepsilon^{\frac1{n+1}}$
for $r\in[a_k,a_{k+1}]$, or
$h(r)>1+n\varepsilon^{\frac1{n+1}}$
for $r\in[a_k,a_{k+1}]$. In any case, using  respectively \eqref{omegaest} and $\omega \le 1$, it follows that
$$
\int_{a_k}^{a_{k+1}}r^{n-1}|h(r)-\omega(r)|\,dr
>\int_{a_k}^{a_{k+1}}r^{n-1}n\varepsilon^{\frac1{n+1}}\,dr
> n^{n+1}\varepsilon,
$$
which from (\ref{omegaint}) contradicts the condition on $h$, and hence proves (\ref{aiaj}).

Since $e^{-2t}<1-t$ and $e^t<1+2t$
for $t\in(0,\frac12)$, (\ref{aiaj}) yields that
$$
e^{-22n\varepsilon^{\frac{1}{n+1}}}\leq 1-11n\varepsilon^{\frac1{n+1}}\leq 
h(a_i),h(a_j)\leq 1+n\varepsilon^{\frac1{n+1}}\le e^{n\varepsilon^{\frac{1}{n+1}}},
$$
 thus $h(a_i)<h(a_j)e^{23(a_j-a_i)}$, and $h(0)<h(a_j)e^{23 a_j}$  by
the log-concavity of $h$. Using the bounds on $h(a_j)$ and $a_j$, we get $h(0)<e^{116n \varepsilon^{\frac1{n+1}} }<
1+250 n\varepsilon^{\frac1{n+1}}.$
On the other hand, the argument leading to (\ref{aiaj})
yields some integer $m\in [1,5]$ such that
$h(a_{-m})\geq 1-11n\varepsilon^{\frac1{n+1}}$. We conclude
by the log-concavity of $h$ that
$$
h(0)\geq \min\{h(a_{-m}),h(a_j)\}\geq 1-11n\varepsilon^{\frac1{n+1}}.
\mbox{ \ \ \proofbox}
$$

\begin{prop}
\label{Banachcentroid}
If the origin $0$ is the centroid of a convex body $K$ in $\R^n$,
and $E\subset K-w\subset (1+\mu)E$ for an $0$-symmetric ellipsoid $E$ and $w\in K$, then
$$
(1-\mu\sqrt{n+1})E\subset K\subset (1+2\mu\sqrt{n+1})E,
$$
holds whenever $\mu\in(0,1/(n+1))$.
\end{prop}
\proof We may assume that $E=B^n$ and $w\neq 0$.
Let $w_0=w/|w|$, and let $B^+$ be the half-ball $\{x\in B^n:\,\langle x,w_0\rangle\geq 0\}$. If $\mu<\frac1{n+1}$, then
$(1+\mu)^{n+1}<e^{\mu(n+1)}<1+2\mu(n+1)$, thus
\begin{eqnarray*}
0&=&\int_{K}\langle x,w \rangle\,dx=
V(K)\langle w,w\rangle+\int_{K-w}\langle x,w\rangle\,dx \\
&>& V(B^n)\langle w,w\rangle+\int_{B^+}\langle x,w\rangle\,dx 
-(1+\mu)^{n+1}\int_{B^+}\langle x,w\rangle\,dx\\
&>&V(B^n)|w|^2-2(n+1)\left(\int_{B^+}\langle x,w_0\rangle\,dx\right)\mu
\cdot |w|.
\end{eqnarray*}
Therefore 
\begin{eqnarray*}
 |w| &\le & (n+1)\mu \frac{\int_{B^n} |\langle x,w_0\rangle|\, dx}{V(B^n)}
      \le  (n+1)\mu \left(\frac{\int_{B^n} \langle x,w_0\rangle^2 \, dx}{V(B^n)}\right)^{\frac12}\\
      &=&     (n+1)\mu \left(\frac{\int_{B^n} |x|^2 \, dx}{nV(B^n)}\right)^{\frac12}\le \mu \sqrt{n+1}.
\end{eqnarray*}
Combining this with our hypothesis $w+B^n\subset K\subset w+(1+\mu)B^n$ readily gives the claim.
 \proofbox

 Now let us prove  Theorem~\ref{FraMeystab}. It is sufficient to consider the case
 $\varepsilon\in(0,\varepsilon_0)$
where $\varepsilon_0>0$ depends on $n$.
Replacing also $f(x)$ by $f(x +z)$ and $g(y)$ by $g(y +z)$
 we may assume that $z=0$.
 For suitable
$\nu,\mu,\lambda>0$, replacing $\varrho(r)$ by
$\nu\varrho(\lambda^2 r)$,
$f(x)$ by $\mu\nu f(\lambda x)$ and $g(x)$ by $(\nu/\mu)g(\lambda x)$,
we may assume that
$$
\int_{\R_+}r^{n-1}\varrho(r^2)\,dr=1\mbox{ \ and \ }
\varrho(0)=f(0)=1.
$$
 Consider the body
$$  K_{f}=\left\{x\in\R^n:\,\int_0^{+\infty} r^{n-1}f(rx)\,dr\geq 1\right\},$$
which is convex since $f$ is log-concave \cite{Bal88}. Its radial function
$$ \|x\|_{K_f}^{-1}=\rho_{K_f}(x):= \sup\big\{t\ge 0; \; tx\in K_f \big\}, \quad x\in S^{n-1}$$
is equal to  $ \left(\int_{\R^+} r^{n-1}f(rx) \, dr \right)^{1/n}.$
Hence, using polar coordinates
shows that
\begin{equation}
\label{Kf}
\int_{\R^n} f(x)\,dx=nV(K_f).
\end{equation}

For $x\in \R^n\backslash \{0\}$, let $f_x,g_x:\R\to\R_+$ be defined by
$f_x(r)=|r|^{n-1}f(rx)$ and $g_x(r)=|r|^{n-1}g(rx)$.
If $\langle x,y\rangle>0$, then the condition on $f,g,\varrho$ yields that
$f_x(r)\cdot g_y(s)\leq m_{xy}(\sqrt{rs})^2$ for
$m_{xy}(r)=r^{n-1}\varrho(r^2\langle x,y\rangle)$
and $r,s\in\R_+$. We deduce by
the Borell-Pr\'ekopa-Leindler inequality Theorem~\ref{PLB} that
$$
\int_{\R_+}f_x(r)\,dr\cdot \int_{\R_+}g_y(r)\,dr \leq
\left(\int_{\R_+}r^{n-1}\varrho(r^2\langle x,y\rangle)\,dr\right)^2
=\langle x,y\rangle^{-n},
$$
and hence
\begin{equation}
\label{KgKfpolar}
K_g\subset K_f^\circ.
\end{equation}
The hypothesis of the theorem  translated in terms of $K_f$ gives
\begin{eqnarray}
\nonumber
 n^2V(B^n)^2&= &
\left(\int_{\R^n}\varrho(|x|^2)\,dx\right)^2
\leq (1+\varepsilon)
\int_{\R^n}f(x)\,dx \int_{\R^n} g(x)\,dx \\
\label{KfKg}
&=&(1+\varepsilon)n^2 V(K_f)V(K_g)\leq (1+\varepsilon) n^2 V(K_f)V(K_f^\circ).
\end{eqnarray}
From the stability version Theorem~\ref{BSstab} of the
Blaschke-Santal\'o inequality, for some $\gamma>0$, $\delta_{\rm BM}(K_f,B^n)<\gamma\,\varepsilon^{\frac1{5n}}$.
 Thus replacing $f(x)$ by
$f(Tx)$ and $g(y)$ by $g(T^{-1}y)$
for a suitable positive definite matrix if necessary,
and applying Proposition~\ref{Banachcentroid}, we may assume that
\begin{equation}
\label{Kfannulus}
B^n\subset K_f\subset (1+O(\varepsilon^{\frac1{5n}}))B^n.
\end{equation}
Using (\ref{KgKfpolar}) we get $K_g\subset K_f^\circ\subset B^n$ and (\ref{KfKg}) yields
\begin{equation}
\label{Kgvolume}
V(K_g)\geq (1+\varepsilon)^{-1}V(B^n)^2V(K_f)^{-1}
\geq (1+O(\varepsilon^{\frac1{5n}}))^{-1}V(B^n).
\end{equation}
For $x\in S^{n-1}$, $\rho_{K_f}(x)=(\int_{\R_+}f_x(r)\,dr)^{\frac1n}\ge 1$
and $\rho_{K_g}(x)=(\int_{\R_+}g_x(r)\,dr)^{\frac1n}\le 1$. We define
\begin{eqnarray*}
\varphi(x)&:=&\int_{\R}f_x(r)\,dr-2=\rho_{K_f}(x)^n+\rho_{K_f}(-x)^n-2 \ge 0, \\
\psi(x)&:=&2-\int_{\R}g_x(r)\,dr=2-\rho_{K_g}(x)^n-\rho_{K_g}(-x)^n \ge 0.
\end{eqnarray*}
In particular (\ref{Kfannulus}) and (\ref{Kgvolume}) yield
\begin{equation}
\label{phipsiint}
\begin{array}{rcl}
\mbox{$\displaystyle \int_{S^{n-1}}$}
\varphi(x)\,dx&=&2n(V(K_f)-V(B^n))=
O(\varepsilon^{\frac1{5n}})\\[1.5ex]
\mbox{$\displaystyle \int_{S^{n-1}}$}
\psi(x)\,dx&=&2n(V(B^n)-V(K_g))=
O(\varepsilon^{\frac1{5n}}),
\end{array}
\end{equation}
where the integration is with respect to the Hausdorff measure on the sphere.
To estimate  $\varphi$ pointwize from above, we use the inclusion (\ref{Kfannulus}).
In order to estimate $\psi$,  we use
(\ref{phipsiint}) and the fact
that a cap of $B^n$ of height $h\le 1$ is of volume larger than
$h^{\frac{n+1}2}V(B^{n-1})/n$ (which forces the  convex subset $K_g$ of the unit ball, with almost the same volume, to 
have a radial function close to $1$ pointwize). More precisely, we obtain that 
there exists $\gamma_0>0$ depending only on $n$, such that
\begin{equation}
\label{phipsimax}
\begin{array}{rcl}
\varphi(x)&\leq&  \gamma_0\varepsilon^{\frac1{5n}}
\mbox{ \ for any $x\in S^{n-1}$},\\[1ex]
\psi(x)&\leq&  \gamma_0\varepsilon^{\frac2{5n(n+1)}}
\mbox{ \ for any $x\in S^{n-1}$}.
\end{array}
\end{equation}
If $\varepsilon_0$ is chosen small enough (depending on $n$),
then (\ref{phipsimax}) yields that
both $\varphi(x)<\frac12$ and $\psi(x)<\frac12$
for any $x\in S^{n-1}$.

Let $x\in S^{n-1}$, and hence
$$
\int_{\R_+}f_x(r)\,dr\geq 1
\mbox{ \ and \ }\int_{\R_+}g_x(r)\,dr\geq
1-\psi(x)\geq (1+2\psi(x))^{-1}.
$$
We define
$m(r)=r^{n-1}\varrho(r^2)$, which
satisfies that $m(e^t)$ is log-concave, and
$f_x(r)\cdot g_x(s)\leq m(\sqrt{rs})^2$ for $r,s\in\R_+$.
Since
$$
\left(\int_{\R_+}m(r)\,dr\right)^2=1
\leq (1+2\psi(x))\int_{\R_+}f_x(r)\,dr\cdot\int_{\R_+}g_x(r)\,dr,
$$
it follows from Theorem~\ref{PLBstab} that there exists
$\alpha(x),\beta(x)>0$ and an absolute constant $c_0>0$
such that
\begin{eqnarray}
\label{fabm}
\int_{\R_+}|\alpha(x)\,f_x(\beta(x)r)-m(r)|\,dr&\leq & c_0\psi(x)^{\frac5{16}}\\
\label{gabm}
\int_{\R_+}|\alpha(x)^{-1}g_x(\beta(x)^{-1}r)-m( r)|\,dr&\leq &
c_0\psi(x)^{\frac5{16}}.
\end{eqnarray}
Using $1\leq \int_{\R_+}f_x(r)\,dr<1+\varphi(x)$ and (\ref{fabm}), we deduce that
\begin{eqnarray*}
\frac{\alpha(x)}{\beta(x)}&\leq& \frac{\alpha(x)}{\beta(x)}\cdot \int_{\R_+}f_x(r)\,dr
=\int_{\R_+}\alpha(x)\,f_x(\beta(x)r)\,dr\\
&\leq& \int_{\R_+}m(r)\,dr+
c_0\psi(x)^{\frac5{16}}=1+c_0\psi(x)^{\frac5{16}}\\
\frac{\alpha(x)}{\beta(x)}&\geq& \frac{\alpha(x)}{\beta(x)}\cdot
(1-\varphi(x))\int_{\R_+}f_x(r)\,dr\\
&=&(1-\varphi(x))\int_{\R_+}\alpha(x)\,f_x(\beta(x)r)\,dr\\
&\geq& (1-\varphi(x))\left(\int_{\R_+}m(r)\,dr-
c_0\psi(x)^{\frac5{16}}\right)\\
&\geq &1-O\left(\max\{\varphi(x),\psi(x)^{\frac5{16}}\}\right).
\end{eqnarray*}
For $a(x)=\alpha(x)^{-1}$ and $b(x)=\beta(x)^{-1}$, we have
\begin{equation}
\label{a(x)/b(x)}
1-c_0\psi(x)^{\frac5{16}}\le \frac{a(x)}{b(x)}\le
1+O\left(\max\{\varphi(x),\psi(x)^{\frac5{16}}\}\right).
\end{equation}
Since $\varphi(x)$ and $\psi(x)$ are even, (\ref{fabm}) can be written in the form
\begin{equation}\label{eq:fabm}
\int_{\R}|r|^{n-1}|f(xr)-a(x)b(x)^{n-1}\,\varrho(b(x)^2 r^2)|\,dr
\le c_0\psi(x)^{\frac5{16}}\frac{a(x)}{b(x)}.
\end{equation}
Thus the hypotheses of Proposition~\ref{logconcave} are satisfied for the log-concave functions $h(r)=f(rx)$ and
$\omega(r)=a(x)b(x)^{n-1}\varrho(b(x)^2 r^2)$ because
$\int_\R |r|^{n-1}\omega(r)dr=\frac{a(x)}{b(x)}$. As $\varrho(0)=f(0)=1$ we get that (using $n+1\leq 2n$)
\begin{equation}
\label{a(x)b(x)}
\left|a(x)b(x)^{n-1}-1\right|=O\left(\psi(x)^{\frac5{32n}}\right).
\end{equation}
We deduce by comparing (\ref{a(x)/b(x)}) and (\ref{a(x)b(x)})
that
$$
\left|a(x)-1\right|=O\left(\max\{\varphi(x),\psi(x)^{\frac5{32n}}\}\right)
\mbox{ \ and \ }
\left|b(x)-1\right|=O\left(\max\{\varphi(x),\psi(x)^{\frac5{32n}}\}\right).
$$
We claim that for any $x\in S^{n-1}$, we have
\begin{eqnarray}
\label{fabm0}
\int_{\R_+}r^{n-1}|f(xr)-\varrho(r^2)|\,dr
&\leq & O\left(\max\{\varphi(x),\psi(x)^{\frac5{32n}}\}\right)\\
\label{gabm0}
\int_{\R_+}r^{n-1}|g(xr)-\varrho(r^2)|\,dr&\leq &
O\left(\max\{\varphi(x),\psi(x)^{\frac5{32n}}\}\right).
\end{eqnarray}
To prove (\ref{fabm0}), we observe 
\begin{eqnarray*}
\int_{\R_+}r^{n-1}|f(xr)-\varrho(r^2)|\,dr
&\le& \int_{\R_+}r^{n-1}|f(xr)-a(x)b(x)^{n-1}\varrho(b(x)^2r^2)|\,dr\\
&+&\int_{\R_+}r^{n-1}a(x)b(x)^{n-1}|\varrho(b(x)^2r^2)-\varrho(r^2)|\,dr\\
&+& \int_{\R_+}r^{n-1}\varrho(r^2)|a(x)b(x)^{n-1}-1|\,dr.
\end{eqnarray*}
Here the first term is $O(\psi(x)^{\frac5{16}})$ by (\ref{eq:fabm}),
and the third term is $O(\psi(x)^{\frac5{32n}})$
by (\ref{a(x)b(x)}). To bound the second term, we first use (\ref{a(x)b(x)}) to get rid of $a(x)b(x)^{n-1}$. To simplify the notations, we put $M=|b(x)^2-1|$. Since $1-M\le b^2\le 1+M$ and $\varrho$ is non-increasing, we obtain 
$$
|\varrho(b(x)^2r^2)-\varrho(r^2)|\le \varrho((1-M)r^2)-\varrho((1+M)r^2).
$$
Thus
\begin{eqnarray*}
\lefteqn{ \int_{\R_+}r^{n-1}|\varrho(b(x)^2r^2)-\varrho(r^2)|\,dr}\\
&\le& \int_{\R_+}r^{n-1}\varrho((1-M)r^2)\,dr-\int_{\R_+}r^{n-1}\varrho((1+M)r^2)\,dr\\
&=&(1-M)^{-\frac{n}{2}}-(1+M)^{-\frac{n}{2}}
=O\left(\max\{\varphi(x),\psi(x)^{\frac5{32n}}\}\right),
\end{eqnarray*}
which in turn yields (\ref{fabm0}). The proof of (\ref{gabm0}) is similar.

Now using  H\"older's inequality and (\ref{phipsiint}), we deduce that
\begin{eqnarray*}
\int_{S^{n-1}}\left(\varphi(x)+\psi(x)^{\frac5{32n}}\right)\,dx&\leq&
 \int_{S^{n-1}} \varphi(x)\,dx +
O\left(\int_{S^{n-1}} \psi(x)\,dx\right)^{\frac5{32n}}\\
&\leq&O(\varepsilon^{\frac1{32n^2}}).
\end{eqnarray*}
Therefore integrating (\ref{fabm0}) and (\ref{gabm0})
over $x\in S^{n-1}$, we have
\begin{eqnarray*}
\int_{\R^n}|f(x)-\varrho(|x|^2)|\,dx
&\leq & O(\varepsilon^{\frac1{32n^2}})\\
\int_{\R^n}|g(x)-\varrho(|x|^2)|\,dx&\leq & O(\varepsilon^{\frac1{32n^2}}).
\end{eqnarray*}
In turn we conclude Theorem~\ref{FraMeystab}.


\section{Proof of Theorem~\ref{Legendrestab} ($\varphi$ convex)}

During the proof of Theorem~\ref{Legendrestab},
 $\gamma_1,\gamma_2,\ldots$ denote
positive constants that depend only on $n$.
We always assume that $\varepsilon\in(0,\varepsilon_0)$,
where $\varepsilon_0>0$ depends on $n$ and $\varrho$,
and is small enough for the argument to work.

We start with some simplification.
We may assume that $\varrho(0)=1=\int_{\R_+}\varrho$. Using the same argument on 
$\varrho$ that the one that we used on $\omega$ in the beginning of the proof of Proposition~\ref{logconcave}, 
 there exists some $t_0\ge 1$ such that $\varrho(t)\leq e^{-t}$
if $t\geq t_0$, and $\varrho(t)\geq e^{-t}$
if $t\in(0,t_0)$. It follows that $\varrho'(0)\geq -1$, and
$$
\int_{\R_+}r^{n-1}\varrho(r^2)\,dr\le \int_{\R_+}r^{n-1}e^{-r^2}\,dr=\frac{\Gamma(n/2)}{2} .
$$

For the log-concave
function $f(x)=\varrho(\varphi(x))$, we may assume that
the origin $0$ is the centre of mass of $K_{f,0}$,
and hence we only check the condition in Theorem~\ref{Legendrestab}
at $z=0$. For $\psi(x)={\cal L}_0\varphi(x)$, let $g(x)=\varrho(\psi(x))$.
It follows from the definition of the Legendre transform that
\begin{equation}
\label{phipsi}
\varphi(x)+\psi(y)\geq \langle x,y\rangle \mbox{ \ for all $x,y\in\R^n$}.
\end{equation}
In particular
$$
f(x)g(y)=\varrho(\varphi(x))\varrho(\psi(y))\leq \varrho^2\left(\frac{\varphi(x)+\psi(y)}2\right)
\leq \varrho^2\left(\frac{\langle x,y\rangle}2\right).
$$
Thus we may apply Theorem~\ref{FraMeystab}, which yields the existence of
 $\xi>0$ and a positive definite matrix
$T:\R^n\to\R^n$, such that
\begin{eqnarray*}
\int_{\R^n}\left|\varrho(|x|^2/2)-\xi\,\varrho(\varphi(Tx))\right|\,dx
&<&\gamma_1\varepsilon^{\frac1{32n^2}} \\
\int_{\R^n}\left|\varrho(|x|^2/2)-\xi^{-1}\varrho(\psi(T^{-1}x))\right|\,dx
&<&\gamma_1\varepsilon^{\frac1{32n^2}},
\end{eqnarray*}
where $\gamma_1$ depends on $n$. Since ${\cal L}_0(\varphi\circ T)=\psi\circ T^{-1}$,
we may assume that $T$ is the identity matrix. 
We choose $R(\varepsilon)$
in a way such that
$$
\varrho(R(\varepsilon)^2)=\varepsilon^{\frac1{64n^2}}.
$$
As $\varrho(t)\leq e^{-t}$ for $t\geq t_0$, it follows 
 that provided  $\varepsilon_0$ is small enough, 
\begin{equation}
\label{Rest}
30<R(\varepsilon)\leq\sqrt{|\log \varepsilon|}/(8n).
\end{equation}
Let $c=-\log \xi$ and $\alpha(x)=-\log \varrho(x)$. Hence $\alpha$
is convex and increasing with $\alpha(0)=0$, $\alpha'(0)\leq 1$, 
where $\alpha'(x)$ denotes the right-derivative. We deduce
\begin{eqnarray*}
\int_{\sqrt{2}R(\varepsilon)B^n}e^{-\alpha(|x|^2/2)}\left|e^{\alpha(|x|^2/2)-\alpha(\varphi(x))-c}-1\right|\,dx
&<&\gamma_1\varepsilon^{\frac1{32n^2}} \\
\int_{\sqrt{2}R(\varepsilon)B^n}e^{-\alpha(|x|^2/2)}\left|e^{\alpha(|x|^2/2)-\alpha(\psi(x))+c}-1\right|\,dx
&<&\gamma_1\varepsilon^{\frac1{32n^2}},
\end{eqnarray*}
which in turn yields by the definition of $R(\varepsilon)$ that
\begin{eqnarray}
\label{etophi}
\int_{\sqrt{2}R(\varepsilon)B^n}
\left|e^{\alpha(|x|^2/2)-\alpha(\varphi(x))-c}-1\right|\,dx
&<&\gamma_1\varepsilon^{\frac1{64n^2}} \\
\label{etopsi}
\int_{\sqrt{2}R(\varepsilon)B^n}
\left|e^{\alpha(|x|^2/2)-\alpha(\psi(x))+c}-1\right|\,dx
&<&\gamma_1\varepsilon^{\frac1{64n^2}}.
\end{eqnarray}

\bigskip
Next we plan to get rid of the exponential function
in (\ref{etophi}) and (\ref{etopsi}).
Define $\tilde{\alpha}(x)=\alpha(|x|^2/2)$.
Then for all $x\in  1.3R(\varepsilon)B^n$,
$$ |\nabla \tilde{\alpha}(x)|=|x|\, \alpha'(|x|^2/2)   \le 1.3 R(\varepsilon) \alpha'\big(0.845 R^2(\varepsilon)\big).$$
Using, for $s,t\ge 0$, the convexity bound $\alpha'(s)\le \frac{\alpha((1+t)s)-\alpha(s)}{ts}\le \frac{\alpha((1+t)s)}{ts}$ together with the relation 
 $\alpha(R(\varepsilon)^2)=|\log\varepsilon|/(64n^2)$, we deduce that 
the function $\tilde{\alpha} $ satisfies
\begin{equation}
\label{tildealpha}
|\nabla \tilde{\alpha}(x)|\leq \gamma_2|\log\varepsilon|
\mbox{ \ \ for $x\in 1.3R(\varepsilon)B^n$}.
\end{equation}

\bigskip
We claim that the convex function $\tilde{\varphi}=\alpha\circ\varphi$ satisfies
\begin{equation}
\label{tildephi}
|\nabla \tilde{\varphi}(x)|\leq 32\gamma_2|\log\varepsilon|
\mbox{ \ \ for $x\in 1.2R(\varepsilon)B^n$}.
\end{equation}
Suppose, to the contrary, that there exists $x_0\in 1.2R(\varepsilon)B^n$
such that the vector $w:=\nabla \tilde{\varphi}(x_0)$ satisfies $|w|> 32\gamma_2|\log\varepsilon|$.
Since $R(\varepsilon)>30$, it follows by (\ref{tildealpha}) that
\begin{equation}
\label{tildealphavariation}
|\tilde{\alpha}(x)-\tilde{\alpha}(x_0)|\leq 3\gamma_2|\log\varepsilon|
\mbox{ \ if $|x-x_0|\leq 3$}.
\end{equation}
We define
$$
\Xi=\Big\{x\in\R^n:\,|x-x_0|\leq 1
\mbox{ \ and } \langle w,x-x_0\rangle\geq \frac12|w|\cdot|x-x_0| \Big\}
\subset 1.3R(\varepsilon)B^n.
$$
If $\tilde{\varphi}(x)\leq\tilde{\alpha}(x_0)-c-4\gamma_2|\log\varepsilon|$
for all $x\in\Xi$, then (\ref{tildealphavariation}) yields
$$
\int_{\sqrt{2}R(\varepsilon)B^n}
\left|e^{\alpha(|x|^2/2)-\alpha(\varphi(x))-c}-1\right|\,dx
>\int_{\Xi}|e^{\gamma_2|\log\varepsilon|}-1|\,dx>
\gamma_1\varepsilon^{\frac1{64n^2}},
$$
provided that $\varepsilon_0$ is small enough. This contradiction to \eqref{etophi}
provides a $y_0\in\Xi$ such that
$\tilde{\varphi}(y_0)\geq\tilde{\alpha}(x_0)-c-4\gamma_2|\log\varepsilon|$.
For $v=\nabla\tilde{\alpha}(y_0)$, we have
\begin{eqnarray*}
\langle v,x_0-y_0\rangle&\leq& \tilde{\varphi}(x_0)-\tilde{\varphi}(y_0)\leq
\langle w,x_0-y_0\rangle \\
&\le &-\frac12 |w|\,|x_0-y_0| \leq -16\gamma_2|\log\varepsilon|\cdot|x_0-y_0|.
\end{eqnarray*}
In particular $|v|\geq 16\gamma_2|\log\varepsilon|$. Next let
$$
\Xi'=\Big\{x\in\R^n:\,1\leq |x-y_0|\leq 2
\mbox{ \ and \ } \langle v,x-y_0\rangle\geq \frac12|w|\cdot|x-y_0| \Big \}
\subset 1.3R(\varepsilon)B^n.
$$
Combining the above definitions and \eqref{tildealphavariation} yields for any $x\in \Xi'$,
 \begin{eqnarray*}
 \tilde{\varphi}(x)&\ge& \tilde\varphi(y_0)+ \langle x-y_0,v\rangle \\
  &\ge & \tilde\alpha(x_0)-c-4\gamma_2 |\log\varepsilon| +\frac12 |w|\,|x-y_0| \\
  &\ge & \tilde\alpha(x_0)-c+4\gamma_2 |\log\varepsilon|\\
  &\ge & \tilde{\alpha}(x)-c+\gamma_2|\log\varepsilon|.
\end{eqnarray*}
Consequently,
$$
\int_{\sqrt{2}R(\varepsilon)B^n}
\left|e^{\alpha(|x|^2/2)-\alpha(\varphi(x))-c}-1\right|\,dx
>\int_{\Xi'}|e^{\gamma_2|\log\varepsilon|}-1|\,dx>\gamma_1\varepsilon^{\frac1{64n^2}},
$$
provided that $\varepsilon_0$ is small enough. This contradicts \eqref{etophi}, hence we may 
 conclude (\ref{tildephi}).

\bigskip
Next we prove that
\begin{equation}
\label{phibig}
\alpha\big(|x|^2/2\big)-\alpha(\varphi(x))-c>-1
\mbox{ \ if $x\in 1.1R(\varepsilon)B^n$}.
\end{equation}
Otherwise suppose that $x_1\in 1.1R(\varepsilon)B^n$
and $\tilde{\alpha}(x_1)-\tilde{\varphi}(x_1)-c\leq-1$.
If $|x-x_1|\leq (96\gamma_2|\log\varepsilon|)^{-1}$ and $\varepsilon_0$ is small enough,
then (\ref{tildealpha}) and (\ref{tildephi}) imply
that $\tilde{\alpha}(x)-\tilde{\varphi}(x)-c\leq-1/3$. Therefore
\begin{eqnarray*}
\int_{\sqrt{2}R(\varepsilon)B^n}
\left|e^{\alpha(|x|^2/2)-\alpha(\varphi(x))-c}-1\right|\,dx
&\ge& \int \left|e^{-\frac13}-1 \right|\, \mathbf1_{|x-x_1|\leq (96\gamma_2|\log\varepsilon|)^{-1}} \, dx\\
&\ge & \gamma_3|\log\varepsilon|^{-n}\, >\, \gamma_1\varepsilon^{\frac1{64n^2}},
\end{eqnarray*}
provided that $\varepsilon_0$ is small enough. This is a contradiction, hence  \eqref{phibig} holds.

\bigskip
Since $|t|<2|e^t-1|$ if $t\geq -1$, combining
(\ref{etophi}) and (\ref{etopsi}) with (\ref{phibig})
and its analogue for $\psi$, we deduce
\begin{eqnarray}
\label{alphaphi}
\int_{R(\varepsilon)B^n}
\left|\alpha(\varphi(x))-\alpha(|x|^2/2)+c\right|\,dx
&<&2\gamma_1\varepsilon^{\frac1{64n^2}} \\
\label{alphapsi}
\int_{R(\varepsilon)B^n}
\left|\alpha(\psi(x))-\alpha(|x|^2/2)-c\right|\,dx
&<&2\gamma_1\varepsilon^{\frac1{64n^2}}.
\end{eqnarray}
For $x\in\R^n$, we define
$
C(x) = \varphi(x)-\frac{|x|^2}2$, $\widetilde{C}(x) = \psi(x)-\frac{|x|^2}2$ and  
\begin{equation}\label{Fpositive}
F(x) =  C(x)+\widetilde{C}(x)\geq 0,
\end{equation}
where the inequality is a consequence of (\ref{phipsi}).
Summing up \eqref{alphaphi} and \eqref{alphapsi}, and using the  convexity of $\alpha$ in the form 
$\alpha(b)-\alpha(a)\ge (b-a)\alpha'(a)$ yields that
\begin{eqnarray*}
4\gamma_1\varepsilon^{\frac1{64n^2}} &\ge &  
  \int_{R(\varepsilon)B^n}
\left(\alpha(\varphi(x))-\alpha(|x|^2/2)+\alpha(\psi(x))-\alpha(|x|^2/2)\right)\,dx \\
&\ge &  \int_{R(\varepsilon)B^n} \alpha'(|x|^2/2)
\left(\varphi(x)-|x|^2/2+\psi(x)-|x|^2/2\right)\,dx \\
&=&  \int_{R(\varepsilon)B^n} \alpha'(|x|^2/2) F(x) \,dx 
\end{eqnarray*}
This is the point where $\alpha$ influences
the estimates. Using \eqref{Fpositive}, we get that 
\begin{equation}
\label{Fest}
\int_{R(\varepsilon)B^n} F(x)\,dx
<\frac{4\gamma_1}{\alpha'(0)}\cdot\varepsilon^{\frac1{64n^2}}.
\end{equation}

\bigskip

Observe that with our notation, \eqref{phipsi} reads as $C(y)+\widetilde C(x)+|x-y|^2/2\ge 0$ or equivalently
$ C(x)\le C(y)+F(x)+|x-y|^2/2$. Since $F$ takes non-negative values, we get that 
  for all $x,y\in\R^n$, 
\begin{equation}\label{eq:C}
|C(x)-C(y)|\leq F(x)+F(y)+\frac{|x-y|^2}{2} \cdot
\end{equation}
 For $t\in \R$, we write $\lceil t \rceil$ for
the smallest integer not smaller than $t$, which satisfies
$\lceil t \rceil+1\leq 2t$ if $\lceil t \rceil\geq 3$.
Set
\begin{equation}\label{defk}
k=\left\lceil\sqrt{\frac{4V(B^n)}{2^{n+1}}}
\left(\int_{R(\varepsilon)B^n}F(z)\,dz \right)^{-\frac12}
\cdot R(\varepsilon)^{\frac{n+2}2}\right\rceil,
\end{equation}
which is at least $3$ if $\varepsilon_0$ is chosen small enough
by (\ref{Rest}) and (\ref{Fest}).
Let us denote
$\sigma:=V(R(\varepsilon)B^n)^{-1}\int_{R(\varepsilon)B^n}C(y)\,dy$.
Taking advantage of \eqref{eq:C}, we get that
\begin{eqnarray}
\nonumber
\int_{R(\varepsilon)B^n}|C(x)-\sigma|\,dx&\leq& V(R(\varepsilon)B^n)^{-1}
\int_{R(\varepsilon)B^n}\int_{R(\varepsilon)B^n}|C(x)-C(y)|\,dxdy\\
\nonumber
&\leq& V(R(\varepsilon)B^n)^{-1}
\sum_{i=1}^k\int_{R(\varepsilon)B^n}\int_{R(\varepsilon)B^n}\\
\nonumber
&&\left|C\left(\mbox{$\frac{i}k\,x+(1-\frac{i}k)y$}\right)-
C\left(\mbox{$\frac{i-1}k\,x+(1-\frac{i-1}k)y$}\right)\right|\,dxdy\\
\nonumber
&\leq&
\sum_{i=0}^k\frac2{V(R(\varepsilon)B^n)}\int_{R(\varepsilon)B^n}\int_{R(\varepsilon)B^n}
F\left(\mbox{$\frac{i}k\,x+(1-\frac{i}k)y$}\right)\,dxdy\\
\label{Cxy}
&&+\frac1{V(R(\varepsilon)B^n)}\sum_{i=1}^k\int_{R(\varepsilon)B^n}\int_{R(\varepsilon)B^n}
\frac{|x-y|^2}{k^2}\,dxdy.
\end{eqnarray}
For $i\in\{0,\ldots,k\}$ in (\ref{Cxy}), we claim that
\begin{equation}
\label{FC}
\int_{R(\varepsilon)B^n}\int_{R(\varepsilon)B^n}
F\left(\mbox{$\frac{i}k\,x+(1-\frac{i}k)y$}\right)\,dxdy
\leq 2^nV(R(\varepsilon)B^n)\int_{R(\varepsilon)B^n}F(z)\,dz.
\end{equation}
If $i\geq k/2$, then  for fixed $y$, using the substitution
$z=\frac{i}k\,x+(1-\frac{i}k)y$, we have
\begin{eqnarray*}
\int_{R(\varepsilon)B^n}\int_{R(\varepsilon)B^n}
F\left(\mbox{$\frac{i}k\,x+(1-\frac{i}k)y$}\right)\,dxdy
&= &\frac{k^n}{i^n}
\int_{R(\varepsilon)B^n}\int_{\frac{i}k\,R(\varepsilon)B^n+(1-\frac{i}k)y}
F(z)\,dzdy\\
&\leq &2^nV(R(\varepsilon)B^n)\int_{R(\varepsilon)B^n}F(z)\,dz.
\end{eqnarray*}
If $i< k/2$, then  for fixed $y$, 
we obtain (\ref{FC}) using the substitution
$z=\frac{i}k\,x+(1-\frac{i}k)y$ for fixed $x$.

In (\ref{Cxy}), we use the rough estimate $|x-y|\leq 2R(\varepsilon)$,
and obtain by (\ref{FC}), \eqref{defk} and $k+1\leq 2k$ that
\begin{eqnarray*}
\int_{R(\varepsilon)B^n}|C(x)-\sigma|\,dx&\leq&
(k+1) 2^{n+1}\int_{R(\varepsilon)B^n}F(z)\,dz+
\frac{4V(B^n)R(\varepsilon)^{n+2}}k\\
&\leq& 3c_0\left(\int_{R(\varepsilon)B^n}F(z)\,dz \right)^{\frac12}
R(\varepsilon)^{\frac{n+2}2}
\end{eqnarray*}
where $c_0>0$ is an absolute constant such that
$\sqrt{2^{n+1}4V(B^n)}<c_0$ for $n\geq 2$.
Since $R(\varepsilon)\leq \sqrt{|\log \varepsilon|}$ by (\ref{Rest}),
we deduce by the definition of $C(x)$ and (\ref{Fest}) that
$$
\int_{R(\varepsilon)B^n}|\mbox{$\varphi(x)-\frac{|x|^2}2$}-\sigma|\,dx<
3c_0\sqrt{\frac{4\gamma_1}{\alpha'(0)}}\cdot\varepsilon^{\frac1{128n^2}}
\cdot |\log \varepsilon|^{\frac{n+2}4},
$$
completing the proof the first part of Theorem~\ref{Legendrestab}.

\section{Proof of  Theorem~\ref{Legendrestab} 
($\varphi$ measurable)}

In this section,  $\eta_1,\eta_2,\ldots$ denote
positive constants that depend only on $n$ and $\varrho$.
Since $\int_{\R^n}\varrho({\cal L}_z\varphi(x))\,dx >0$ for all $z$, the function ${\cal L}_z\varphi$
cannot be identically infinite. Hence
we may consider the lower convex hull 
$\varphi_*={\cal L}_z{\cal L}_z\varphi$ of $\varphi$.
It follows that ${\cal L}_z\varphi_*={\cal L}_z\varphi$.
We may assume as in the proof of Theorem~\ref{Legendrestab} that
$\varrho(0)=1=\int_{\R_+}\varrho$, and hence $\varrho'(0)\geq -1$.
Let again $\alpha(t)=-\log\varrho(t)$, which is convex, increasing,
and satisfies $\alpha(0)=0$ and $0<\alpha'(0)\leq 1$,
where $\alpha'(x)$ denotes the right-derivative.

For $t\in\R$, we also introduce 
$$
\alpha_*(t)=\left\{
\begin{array}{rl}
\alpha(t)&\mbox{ \ if $t\geq 0$}\\
\alpha'(0)\cdot t & \mbox{ \ if $t\leq 0$}.
\end{array}\right.
$$
As we shall see shortly, we can replace $\alpha$ by $\alpha_*$ in the
inequalities. Observe first that  $\alpha_*\le \alpha$ and that  
$$
\alpha'_*(t)\geq \alpha'(0)=\alpha'_*(0) 
\mbox{ \ \ for all $t\in\R$}.
$$
Let $\varrho_*(t)=e^{-\alpha_*(t)}$. 
As $\varrho_*(t)\geq \varrho(t)$, $\varphi_*(x)\leq \varphi(x)$
and ${\cal L}_z\varphi_*={\cal L}_z\varphi$, we have
\begin{eqnarray}
\nonumber
\int_{\R^n}\varrho_*(\varphi_*(x))\,dx \int_{\R^n}
\varrho_*({\cal L}_z\varphi_*(x))\,dx&\geq &
\int_{\R^n}\varrho_*(\varphi(x))\,dx \int_{\R^n}
\varrho_*({\cal L}_z\varphi(x))\,dx\\
\nonumber
&\geq & (1-\varepsilon)
\left(\int_{\R^n}\varrho(|x|^2/2)\,dx\right)^2\\
\label{starlow}
&= & (1-\varepsilon)
\left(\int_{\R^n}\varrho_*(|x|^2/2)\,dx\right)^2
\end{eqnarray}
for any $z$.
We may assume that  
the origin $0$ is the centre of mass of $K_{f,0}$
for the log-concave
function $f=\varrho_*\circ \varphi_*$. Therefore
\begin{equation}
\label{starup}
\int_{\R^n}\varrho_*(\varphi_*(x))\,dx 
\int_{\R^n}\varrho_*({\cal L}_0\varphi_*(x))\,dx
\leq \left(\int_{\R^n}\varrho_*(|x|^2/2)\,dx\right)^2.
\end{equation}
We have proved in the course of the argument for 
Theorem~\ref{Legendrestab} that
possibly after a positive definite linear transformation, there
exists $\sigma\in\R$ such that
\begin{equation}
\label{close-quadratic0}
\int_{R_*(\varepsilon)B^n}
\left|\varphi_*(x)-\sigma-\frac{|x|^2}2\right|\,dx
<\eta_1\varepsilon^{\frac1{129n^2}}
\end{equation}
where 
\begin{equation}
\label{alphaR}
\alpha_*(R_*(\varepsilon)^2)=\frac{|\log\varepsilon|}{64n^2}.
\end{equation}
In particular $\lim_{\varepsilon\to 0}R_*(\varepsilon)=+\infty$ and
$30< R_*(\varepsilon)\leq\sqrt{|\log\varepsilon|}/(8n)$.
Set $R(\varepsilon):=\frac12\,R_*(\varepsilon)$.

\begin{prop}

\label{close-quadratic}
If $\varepsilon>0$ is small enough, then
$$
\left|\varphi_*(x)-\sigma-\frac{|x|^2}2\right|
<\eta_2\varepsilon^{\frac2{129n^2(n+2)}}<1
\mbox{ \ \ \ for all $x\in \frac53\,R(\varepsilon)\,B^n$.}
$$
\end{prop}

\proof
Let us denote   by $c$ the convex function $\varphi_*-\sigma$ and $f(x)=c(x)-|x|^2/2$. Set $\delta=\eta_1\varepsilon^{1/(129 n^2)}$.
Assume that $\varepsilon$ is small enough so that $\delta<1$. Our starting point is \eqref{close-quadratic0}
which reads as $\int_{R_*(\varepsilon) B^n} \big|f| \le \delta$. 
 
 Let $r\in (0,1)$ and $x\in \R^n$ with $|x|\le R_*(\varepsilon)-1$.
 If $v(x)$ is a subgradient of $c$ at $x$, we get by convexity of $c$ that for all $y$,
 $$ f(y)\ge f(x)+\langle v(x)-x,y-x\rangle -\frac{|y-x|^2}{2}\cdot$$
 Since the ball $B(x,r)$ of center $x$ and radius $r$ is included in $B(0,R_*(\varepsilon))=R_*(\varepsilon)\, B^n$,
 we deduce that 
 \begin{eqnarray*}
 \delta &\ge & \int_{B(x,r)} |f| \ge \int _{B(x,r)} f(y) \, dy \\
 &\ge & f(x)\, V(B(x,r)) - \int_{B(x,r)} \frac{|y-x|^2}2\, dy =v_n r^n f(x)- c_n r^{n+2}
 \end{eqnarray*}
 for suitable quantities $v_n,c_n$ depending only on $n$.
 Rearranging,  $f(x)\le (\delta+ c_n r^{n+2})/(v_n r^n)$. Choosing $r=\delta^{1/(n+2)}<1$, we obtain that
 for all $x$ with $|x|\le R_*(\varepsilon)-1$,
\begin{equation}\label{eq:upper-f}
 f(x)\le d_n \delta^{\frac{2}{n+2}}.
\end{equation}

 In order to establish the proposition, it remains to prove a similar lower bound on $f$.
Consider $x\in \R^n$ with $|x|\le \R_*(\varepsilon)-2$. Let $r\in(0,1)$ to be specified later. 
Consider a point $y\in B(x,r)\setminus\{x\}$. It can be written as $y=x+su$ with $s \in (0,r]$ and
$u\in\R^n$, $|u|=1$. By convexity,
\begin{eqnarray*}
c(y) &\le & \frac{s}{r} c(x+ru)+\left(1-\frac{s}{r}\right)c(x) \\
   &=&  \frac{s}{r} \left(f(x+ru)+\frac{|x+ru|^2}2\right)+\left(1-\frac{s}{r}\right)\left( f(x)+\frac{|x|^2}2\right). 
\end{eqnarray*}
 Rearranging the squares and using the upper bound \eqref{eq:upper-f} gives
 \begin{eqnarray*}
-f(y) &=  & \frac{|y|^2}2-c(y) \ge -\frac 12 s(r-s)-\frac{s}{r} f(x+ru)-\left(1-\frac{s}{r}\right)f(x) \\
   &\ge &  -\frac 12 s(r-s)-\frac{s}{r}  d_n \delta^{\frac{2}{n+2}}-\left(1-\frac{s}{r}\right)f(x).
\end{eqnarray*}
Integrating in $y=x+su$ in spherical coordinates of origin $x$, and changing variables $s=rt$, $t\in(0,1]$ gives
for suitable positive numbers depending only on the dimension
\begin{eqnarray*} 
 \delta &\ge& \int_{R_*(\varepsilon) B^n} |f| \ge \int_{B(x,r)} |f|\ge \int_{B(x,r)} -f(y)\, dy \\
 &\ge & -\frac 12 \int_0^r s(r-s) nV(B^n)s^{n-1}ds -  d_n \delta^{\frac{2}{n+2}}\int_0^r \frac{s}{r} nV(B^n)s^{n-1}ds\\ && -f(x)\int_0^r \left(1-\frac{s}{r}\right) nV(B^n)s^{n-1}ds \\
 &=& -c_n r^{n+2}-d'_n \delta^{\frac{2}{n+2}} r^n-c'_n f(x) r^n.
\end{eqnarray*} 
Choosing $r=\delta^{1/(n+2)}<1$ and rearranging yields $f(x)\ge -c''_n \delta^{\frac{2}{n+2}}$, provided
$|x|\le R_*(\varepsilon)-2$. Since $R_*(\varepsilon)>30$, the claim follows.
\proofbox

We now estimate how close the weight function $\varrho_*\circ \varphi_*$
is  to be a constant function on $R(\varepsilon)B^n$.
We claim that
\begin{equation}
\label{diamclaim}
\frac{\varrho_*\circ \varphi_*(x)}{\varrho_*\circ \varphi_*(y)}
\geq\varepsilon^{\frac1{16n^2}} \mbox{ \ for all $x,y\in R(\varepsilon)B^n$}.
\end{equation}
Since the function $\alpha_*=-\log \rho_*$ is increasing, we deduce from 
Proposition~\ref{close-quadratic} that  
$$
|\alpha_*(\varphi(x))-\alpha_*(\varphi(y))|\leq
\alpha_*\left(\frac{R(\varepsilon)^2}2+1+\sigma\right)-
\alpha_*(\sigma-1) \mbox{ \ for $x,y\in R(\varepsilon)B^n$}.
$$
Therefore it is sufficient to prove that
\begin{equation}
\label{diamclaim0}
\Omega:=\alpha_*\left(\frac{R(\varepsilon)^2}2+1+\sigma\right)-
\alpha_*(\sigma-1)\leq \frac{|\log\varepsilon|}{16n^2}.
\end{equation}
It follows by (\ref{alphaphi}) and (\ref{alphapsi}) that
$$
\int_{R_*(\varepsilon)B^n}
\Big(\alpha_*\big(\varphi_*(x)\big)+\alpha_*\big({\cal L}_0\varphi_*(x)\big)-2\alpha_*\big(|x|^2/2\big)\Big)\,dx
<\eta_8\varepsilon^{\frac1{64n^2}}. 
$$
We note that by definition $\varphi_*(x)+{\cal L}_0\varphi_*(x)\geq |x|^2$
for all $x\in\R^n$. Hence the monotonicity of $\alpha_*$ yields
\begin{equation}\label{eq:alphaeta8}
\int_{R_*(\varepsilon)B^n}
\Big(\alpha_*\big(\varphi_*(x)\big)+\alpha_*\big(|x|^2-\varphi_*(x)\big)-2\alpha_*\big(|x|^2/2\big)\Big)\,dx
<\eta_8\varepsilon^{\frac1{64n^2}}. 
\end{equation}
Next, we bound from below the three terms appearing inside the above integral, when the variable is in the 
smaller domain $\frac53\,R(\varepsilon)B^n\backslash (\frac43\,R(\varepsilon)B^n)$.
Observe that if 
$x\in \frac53\,R(\varepsilon)B^n\backslash (\frac43\,R(\varepsilon)B^n)$,
then  by 
Proposition~\ref{close-quadratic}
$$\varphi_*(x)\geq\frac{(\frac43\,R(\varepsilon))^2}2+\sigma-1\geq
\frac{R(\varepsilon)^2}2+\sigma+1.$$
Since $\alpha_*$ is convex, increasing and verifies $\alpha_*(0)=0$, $\alpha'_*(0)=\alpha'(0)$ we get that
$$ \alpha_*(\varphi_*(x))\ge \alpha_*\left(\frac{R(\varepsilon)^2}2+\sigma+1\right)
=\Omega + \alpha_*(\sigma-1) \ge \Omega+\alpha'(0)(\sigma-1).$$
Still assuming that $x\in \frac53\,R(\varepsilon)B^n\backslash (\frac43\,R(\varepsilon)B^n)$ and taking advantage
of  Proposition~\ref{close-quadratic}, we obtain that 
$$\alpha_*(|x|^2-\varphi_*(x)) \ge \alpha_*(-\sigma-1)\geq \alpha'(0)(-\sigma-1).$$
Eventually, Equation~(\ref{alphaR}) gives for $x\in \frac53\,R(\varepsilon)B^n\backslash (\frac43\,R(\varepsilon)B^n)$,
$$ \alpha_*(|x|^2/2)\le \alpha_*\big((5R(\varepsilon)/3)^2/2\big)\le  \alpha_*(R_*(\varepsilon)^2)= \frac{|\log\varepsilon|}{64n^2}.$$ 
Since the integrand in \eqref{eq:alphaeta8} is always non-negative, the above three inequalities together
with \eqref{eq:alphaeta8} easily yield that
$$
\int_{\frac53\,R(\varepsilon)B^n\backslash (\frac43\,R(\varepsilon)B^n)}
\left(\Omega-2\alpha'(0)-\frac{|\log\varepsilon|}{32n^2}\right)\,dx
<\eta_8\varepsilon^{\frac1{64n^2}}. 
$$
From this, we  conclude that (\ref{diamclaim0})
and thus (\ref{diamclaim}) hold  if $\varepsilon$ is small enough.

\bigskip
Next, we define the set 
$$
\Psi=: \left\{x\in R(\varepsilon)\,B^n:\,
\varphi(x)>\varphi_*(x)+\varepsilon^{\frac1{128n^2}}
\right\}.
$$ 
Since the inequality $\varphi(x)\geq\varphi_*(x)$ is true for all $x\in\R^n$, 
it follows from equation   (\ref{close-quadratic0}) and the bound $R(\varepsilon)<\sqrt{|\log\varepsilon|}$
 that
$$
\int_{R(\varepsilon)B^n\backslash \Psi}
\left|\frac{|x|^2}2+\sigma-\varphi(x)\right|\,dx
<\eta_{10}\varepsilon^{\frac1{129n^2}}.
$$
Therefore our final task is to provide
a suitable upper bound on the volume of the set $\Psi$.

Let $R_0>0$ be defined by $\alpha'(0)\cdot(\frac{R_0^2}2-2)=1$. Since $\alpha'(0)\in (0,1]$, $R_0\ge \sqrt 6$.
From now on we consider $R\in[R_0,R(\varepsilon)]$.
Since $\alpha'_*(t)\geq \alpha'(0)$ for $t\in\R$, we have for all $x\in\Psi$
\begin{eqnarray*}
\varrho_*(\varphi(x))&=&e^{-\alpha_*(\varphi(x))} \le 
e^{-\alpha_*\left(\varphi_*(x)+\varepsilon^{\frac1{128n^2}}\right)}\\ &\leq&
 e^{-\alpha'(0)\varepsilon^{\frac1{128n^2}}}
\varrho_*(\varphi_*(x))
\le \left(1-\alpha'(0)\varepsilon^{\frac1{127n^2}}\right) \varrho_*(\varphi_*(x)),
\end{eqnarray*}
where the last inequality is valid if $\varepsilon$ is small enough.
This improves on the trivial estimate $\varrho_*\circ \varphi \le \varrho_*\circ \varphi_*$. Let us see how the improvement passes to integrals:
\begin{eqnarray*}
\lefteqn{\int_{R\,B^n}\varrho_*\circ \varphi =  \int_{\Psi \cap R\,B^n}\varrho_*\circ \varphi +
                                           \int_{  R\,B^n\setminus \Psi}\varrho_*\circ \varphi }\\
   &\le &  \left(1-\alpha'(0)\varepsilon^{\frac1{127n^2}}\right) \int_{\Psi \cap R\,B^n}\varrho_*\circ \varphi_* +
                                           \int_{  R\,B^n\setminus \Psi}\varrho_*\circ \varphi_*\\
   &=&     \int_{R\,B^n}\varrho_*\circ \varphi_*  -   \alpha'(0)\varepsilon^{\frac1{127n^2}}   
     \int_{\Psi \cap R\,B^n}\varrho_*\circ \varphi_*.                        
\end{eqnarray*}
However,   (\ref{diamclaim}) readily gives 
$$ \int_{\Psi \cap R\,B^n}\varrho_*\circ \varphi_* \ge \varepsilon^{\frac{1}{16n^2}} \frac{V(\Psi\cap R\, B^n)}{V(R\, B^n)} \int_{R \, B^n} \varrho_*\circ \varphi_* .$$
Hence, combining this with the former estimate, we deduce that 
\begin{eqnarray*}
\int_{R\,B^n}\varrho_*\circ \varphi &\leq &
\left(1-\frac{\varepsilon^{\frac1{16n^2}}V(\Psi\cap R \, B^n)}{V(R\,B^n)}
\cdot \alpha'(0)\varepsilon^{\frac1{127n^2}}\right)
\int_{R\,B^n}\varrho_*\circ \varphi_* .
\end{eqnarray*}
Our goal is to draw information on the volume of $\Psi$ from the above inequality and the almost equality
in the functional Blaschke-Santal\'o inequality. But this requires  a similar inequality
involving integrals on the whole space. Building on the latter estimate,
\begin{eqnarray}\label{eq:total-ball}
\int_{\mathbb R^n}\varrho_*\circ \varphi &=&\int_{R\,B^n}\varrho_*\circ \varphi +\int_{\mathbb R^n \setminus R\,B^n}\varrho_*\circ \varphi \nonumber\\
&\le & 
\left(1-\varepsilon^{\frac1{8n^2}}\alpha'(0)\frac{V(\Psi\cap R \, B^n)}{V(R\,B^n)}
 \right)\int_{R\,B^n}\varrho_*\circ \varphi_*+\int_{\mathbb R^n \setminus R\,B^n}\varrho_*\circ \varphi_* \nonumber\\
  &=& \int_{\mathbb R^n}\varrho_*\circ \varphi_* -\varepsilon^{\frac1{8n^2}}\alpha'(0)\frac{V(\Psi\cap R \, B^n)}{V(R\,B^n)} \int_{R\,B^n}\varrho_*\circ \varphi_*.
\end{eqnarray}
If $|x|=R_0$, then 
Proposition~\ref{close-quadratic} and  the properties of $\alpha_*$    yield
\begin{eqnarray*}
&&\alpha_*(\varphi_*(x))-\alpha_*(\varphi_*(0))\ge
\alpha_*\left(\frac{R_0^2}2+\sigma-1\right)-
\alpha_*(\sigma+1)\\
&&\ge \alpha'_*(\sigma+1) \left( \frac{R_0^2}2-2\right)= \frac{\alpha'_*(\sigma+1)}{\alpha'(0)}\geq 1,
\end{eqnarray*}
thus the log-concave function $\varrho_*\circ\varphi_*$ verifies $\varrho_*(\varphi_*(x))\leq e^{-1}\varrho_*(\varphi_*(0))$
whenever $|x|=R_0$. Then, elementary estimates for one-dimensional log-concave functions (applied on all radii)
  give
$$
\frac{\int_{R_0\,B^n}\varrho_*\circ \varphi_* }{\int_{\R^n}\varrho_*\circ \varphi_*} \geq
\frac{\int_{R_0\,B^n}e^{-|x|/R_0}\,dx}{\int_{\R^n}e^{-|x|/R_0}\,dx}\cdot
$$
Since the latter ratio depends only on $n$, we consider it as a constant. Hence we deduce 
from \eqref{eq:total-ball} that for $R\in [R_0,R(\varepsilon)]$
$$
\frac{\int_{\R^n}\varrho_*\circ \varphi}{\int_{\R^n}\varrho_*\circ \varphi_*}\le
 1-\frac{\eta_{11}\varepsilon^{\frac1{8n^2}}V(\Psi\cap R \, B^n)}{V(R\, B^n)} 
 \cdot
$$
On the other hand,
(\ref{starlow}) and (\ref{starup}) give $\frac{\int_{\R^n}\varrho_*\circ \varphi}{\int_{\R^n}\varrho_*\circ \varphi_*}\ge 1-\varepsilon$. Comparing the latter two estimates leads to 
$$
V(\Psi\cap R\, B^n)\leq \eta_{11}^{-1}\varepsilon^{1-\frac1{8n^2}}V(R\, B^n).
$$
The proof of   Theorem~\ref{Legendrestab} is therefore complete.

\medskip \noindent
F. Barthe : Institut de math\'ematiques de Toulouse, UMR 5219. Universit\'e Paul Sabatier. 31062 Toulouse cedex 9.  FRANCE.

\noindent E-mail: barthe@math.univ-toulouse.fr

\medskip \noindent
K. B\"or\"oczky :  Alfr\'ed R\'enyi Institute of Mathematics,
	 Hungarian Academy of Sciences
	 1053 Budapest, Re\'altanoda u. 13-15.
	 HUNGARY.

\noindent E-mail: carlos@renyi.hu

\medskip \noindent
M. Fradelizi : Universit\'e Paris-Est Marne-la-Vall\'ee.
Laboratoire d'Analyse et de Math\'ematiques Appliqu\'ees UMR 8050.
5 Bd Descartes.
Champs-sur-Marne.
77454 Marne-la-Vall\'ee Cedex 2. FRANCE.

\noindent E-mail: matthieu.fradelizi@univ-mlv.fr

\begin{thebibliography}{9}


\bibitem{AKM04}
S. Artstein-Avidan, B. Klartag, V.D. Milman:
On the Santal\'o point of a function and a
functional Santal\'o inequality,
Mathematika 54 (2004), 33--48.

\bibitem{Bal}
K.M. Ball:
Isometric embedding in $L_p$ and sections of convex sets.
PhD thesis, University of Cambridge, 1988.

\bibitem{Bal88}
K.M. Ball:
Logarithmically concave functions and sections of convex sets in $\R^n$.
Studia Math.,  88  (1988),   69--84.



\bibitem{BB1}
K.M. Ball, K.J. B\"or\"oczky:
Stability of the Pr\'ekopa-Leindler inequality.
Mathematika, accepted.

\bibitem{BB2}
K.M. Ball, K.J. B\"or\"oczky:
Stability of some versions of the Pr\'ekopa-Leindler inequality.
Monatsh. Math., accepted.


\bibitem{Bla17}
W. Blaschke:
\"Uber affine Geometrie VII.
Neue Extremeigenschaften von Ellipse und Ellipsoid.
Leipz. Ber., 69 (1917), 306--318.

\bibitem{Bla22}
W. Blaschke:
\"Uber affine Geometrie  XXXVII.
Eine Versch\"arfung von Minkowskis Ungleichheit
f\"ur den gemischten Fl\"acheninhalt.
Hamb. Abh., 1 (1922), 206--209.



\bibitem{BoF87}
T. Bonnesen, W. Fenchel:
Theory of convex bodies.
BCS Associates,  1987.



\bibitem{Bor75}
C. Borell:
Convex set functions in $d$-space.
Period. Math. Hungar. 6 (1975), 111--136.



\bibitem{Bor}
K.J. B\"or\"oczky:
Stability of the Blaschke-Santal\'o
and the affine isoperimetric inequality.
Adv. Math., 225 (2010), 1914-1928.



\bibitem{FrM07}
M. Fradelizi, M. Meyer:
Some functional forms of Blaschke-Santal\'o inequality.
Math. Z., 256 (2007), 379--395.


\bibitem{Gru07}
P.M. Gruber:
Convex and discrete geometry.
Springer, Berlin, 2007.


\bibitem{Hug96}
D. Hug:
Contributions to affine surface area.
Manuscripta Math., 91 (1996), 283--301.


\bibitem{Leh1}
J. Lehec:  
A direct proof of the functional Santal\'o inequality.  
C. R. Math. Acad. Sci. Paris  347  (2009),  no. 1-2, 55--58.

\bibitem{Leh2}
J. Lehec: 
Partitions and functional Santal\'o inequalities.  
Arch. Math. (Basel)  92  (2009),  no. 1, 89--94..



\bibitem{Lut91}
E. Lutwak:
Extended affine surface area.
Adv. Math, 85 (1991), 39--68.

\bibitem{Lut93}
E. Lutwak:
 Selected affine isoperimetric inequalities.
In: Handbook of convex geometry,
North-Holland, Amsterdam, 1993, 151--176.


\bibitem{MeP90}
M. Meyer, A. Pajor:
On the Blaschke-Santal\'o inequality.
Arch. Math. (Basel) 55 (1990), 82--93.


\bibitem{MeR06}
M. Meyer, S. Reisner:
Shadow systems and volumes of polar convex bodies.
Mathematika, 53 (2006), 129--148.


\bibitem{Pet85}
C.M. Petty:
Affine isoperimetric problems.  Discrete geometry and
convexity (New York, 1982),  113--127,
Ann. New York Acad. Sci., 440, New York Acad. Sci., New York, 1985.



\bibitem{Sai81}
J. Saint-Raymond:
Sur le volume des corps convexes sym\'etriques.
 Initiation Seminar on Analysis: G. Choquet-M. Rogalski-J. Saint-Raymond,
20th Year: 1980/1981, Exp. No. 11, 25 pp.,
Publ. Math. Univ. Pierre et Marie Curie, 46, Univ. Paris VI, Paris, 1981.

\bibitem{San49}
L.A. Santal\'o:
An affine invariant for convex bodies of $n$-dimensional space. (Spanish)
Portugaliae Math., 8 (1949), 155--161.



\bibitem{Sch93}
R. Schneider:
Convex Bodies: The Brunn-Minkowski Theory.
Cambridge University Press, 1993.



\end{thebibliography}
\end{document}